\documentclass[a4paper,12pt]{article}
\usepackage{geometry,latexsym,amssymb,amsmath,amsthm,color,bm}
\usepackage[latin5]{inputenc}
\usepackage{enumerate}
\usepackage{enumitem}
\usepackage[T1]{fontenc}
\usepackage{authblk}
\usepackage[all]{xy}
\usepackage{palatino}
\usepackage{indentfirst}
\usepackage{titlesec}
\usepackage{graphics}

\usepackage{setspace}
\theoremstyle{plain}
\newtheorem{example}{Example}[section]
\newtheorem{Def}[example]{Definition}
\newtheorem{Exam}[example]{Example}
\newtheorem{Prop}[example]{Proposition}
\newtheorem{Theo}[example]{Theorem}
\newtheorem{Lem}[example]{Lemma}
\newtheorem{Rem}[example]{Remark}

\newenvironment{Prf}{{\bf Proof:} } {\hfill $\Box$\mbox{}
\mbox{}}

\def\Lbnz{\bm{\mathsf{Lbnz}}}
\def\LX{\bm{\mathsf{XMod(Lbnz)}}}
\def\bmC{\bm{\mathsf{C}}}
\def\Cat{\bm{\mathsf{Cat}}}
\def\Gpd{\bm{\mathsf{Gpd}}}

\def\LGdC/G{\bm{\mathsf{Cov_{\Cat(\Lbnz)}/G}}}
\def\LGdA(G){\bm{\mathsf{Act_{\Cat(\Lbnz)}(G)}}}
\def\LXMC(L){\bm{\mathsf{Cov_{\XMod(\Lbnz)}/(L_1,L_0,\partial)}}}

\def\wtilde{\widetilde}

\def\St{\mathsf{St}}

\def\Ob{\operatorname{Ob}}

\def\XMod{\bm{\mathsf{XMod}}}

\def\NSGGd/G{\mathsf{NSGGd/G}}
\def\NSCM/(A,B,\alpha){\mathsf{NSCM/(A,B,\alpha)}}
\def\SpPXM/(A,B,\alpha){\mathsf{SpPXM/(A,B,\alpha)}}
\def\LXM/(A,B,\alpha){\mathsf{LXMod/(A,B,\alpha)}}
\def\LXM{\mathsf{LXMod}}

\def\GGdA(G){\bm{\mathsf{GpGpdAct(G)}}}
\def\GGdC/G{\bm{\mathsf{GpGpdCov/G}}}
\def\GGdCov/X{\bm{\mathsf{GpGpdCov/\pi X}}}
\def\GdC/G{\bm{\mathsf{GpdCov/G}}}
\def\TGrC/X{\mathsf{TGrCov/X}}
\def\GdA(G){\bm{\mathsf{GpdAct(G)}}}
\def\Act(G){\bm{\mathsf{GpdAct(G)}}}
\def\Cov/G{\bm{\mathsf{GpdCov/G}}}
\def\C{\bm{\mathsf{C}}}

\def\epsilon{\varepsilon}

\begin{document}
\title{Actions of internal groupoids in the category of Leibniz algebras\thanks{This work has been supported by Research Fund of the Aksaray University. Project Number:2016-019.}}

\author[a]{Tunçar ŞAHAN\thanks{\textbf{Correspondence :} Tunçar ŞAHAN (e-mail : tuncarsahan@gmail.com)}}
\author[b]{Ayhan ERCİYES\thanks{Ayhan ERCİYES (e-mail : ayhan.erciyes@hotmail.com)}}
\affil[a]{\small{Department of Mathematics, Aksaray University, Aksaray, TURKEY}}
\affil[b]{\small{Department of Mathematics and Science Education, Aksaray University, Aksaray, TURKEY}}

\date{}

\maketitle

\begin{abstract}
The aim of this paper is to characterize the notion of internal category (groupoid) in the category of Leibniz algebras and investigate the properties of well-known notions such as covering groupoid and groupoid operations (actions) in this category. Further, for a fixed internal groupoid $G$, we prove that the category of covering groupoids of $G$ and the category of internal groupoid actions of $G$ on Leibniz algebras are equivalent. Finally we interpret the corresponding notion of covering groupoids in the category of crossed modules of Leibniz algebras.
\end{abstract}

\\ {\bf Classification:} 17A32, 20L05, 18D35

\section{Introduction}
Covering groupoids have an important role in the applications of groupoids (see for example  \cite{Br1} and \cite{Hi}). It is well known that for a groupoid $G$, the category $\Act(G)$ of  groupoid actions  of $G$ on  sets, these are also called operations or $G$-sets, are equivalent to  the category $\Cov/G$ of  covering groupoids of $G$. For the topological version of this equivalence, see \cite[Theorem 2]{Br-Da-Ha}.

If $G$ is a group-groupoid, which is an internal groupoid in the category of groups, then the category $\GGdC/G$ of group-groupoid coverings of $G$ is equivalent to the category $\GGdA(G)$ of group-groupoid actions of $G$ on groups \cite[Proposition 3.1]{Br-Mu1}. In \cite{Ak-Al-Mu-Sa} this result has recently generalized to the case where $G$ is an internal groupoid for an algebraic category $\C$, acting on a group with operations. Covering groupoids of a categorical group have been studied in \cite{Mu-Sa}. 

In \cite{BS1} it was proved that the categories of crossed modules and group-groupoids, under the name of $\mathcal{G}$-groupoids, are equivalent (see also  \cite{Loday82} for an alternative equivalence in terms of an algebraic object called {\em cat$^n$-groups}). By applying this equivalence of the categories,  normal and quotient objects in the category of group-groupoids  have been recently obtained in \cite{Mu-Sa-Al}. The study of internal category theory was continued in the works of Datuashvili \cite{Kanex} and \cite{Wh}. Moreover, she developed cohomology theory of internal categories in categories of groups with operations \cite{Dat} and \cite{Coh} (see also \cite{Por} for more information on internal categories in categories of groups with operations). The  equivalences of the categories  in \cite{BS1} enable us to generalize some results on group-groupoids to the more general internal groupoids for a certain  algebraic category $\C$ (see for example \cite{Ak-Al-Mu-Sa}, \cite{Mu-Be-Tu-Na},  \cite{Mu-Tu}  and \cite{Mu-Ak}).

In the mid-nineteenth century, Whitehead introduced the notion of crossed module, in a series of papers \cite{Wth1,Wth3,Wth2}, as algebraic models for (connected) homotopy 2-types (i.e. connected spaces with no homotopy group in degrees above 2), in much the same way that groups are algebraic models for homotopy 1-types. A crossed module consists of groups $A$ and $B$, where $B$ acts on $A$ by automorphisms, and a homomorphism of groups $\alpha\colon A\rightarrow B$ satisfying \textbf{\textit{(i)}} $\alpha( ^{b} a)=b+\alpha(a)-b$ and \textbf{\textit{(ii)}} $ ^{\alpha(a)} a_1 = a+a_1-a$ for all $a,a_1 \in A$ and $b\in B$. Crossed modules can be viewed as 2-dimensional groups \cite{BrLowDim} and have been widely used in: homotopy theory, \cite{BHS}; the theory of identities among relations for group presentations, \cite{BrownHubesshman}; algebraic K-theory \cite{Loday}; and homological algebra, \cite{Hub,Lue}. See \cite[pp.49]{BHS} for some discussion of the relation of crossed modules to crossed squares and so to homotopy 3-types. The equivalence between crossed modules and group groupoids, proved in \cite{BS1} and has been found important in applications. It is generalised in \cite{Por}.

In this paper, first we defined and investigated some properties of internal categories (and hence internal groupoids) in the category of Leibniz algebras. Further we defined coverings and actions in the category of internal groupoids in the category of Leibniz algebras and proved that the category of internal groupoid actions and the category of covering groupoids of a fixed internal groupoid $G$ in the category of Leibniz algebras are equivalent. Finally, using the equivalence of the categories internal groupoids in the category of Leibniz algebras and crossed modules of Leibniz algebras, we interpreted the notion of covering in the category of crossed modules of Leibniz algebras. 

\section{Preliminaries}

A Leibniz algebra $L$ is a $\Bbbk$-vector space equipped with a bilinear map $\left[-,-\right]\colon L\times L\rightarrow L$, satisfying the Leibniz identity $\left[x,\left[y,z\right]\right]=\left[\left[x,y\right],z\right]-\left[\left[x,z\right],y\right]$ for all $x,y,z\in L$. Leibniz algebras are the generalization of Lie algebras. Indeed, for a Leibniz algebra $L$, if $\left[x,x\right]=0$ for all $x\in L$, then $L$ becomes a Lie algebra. On the other hand, every Lie algebra is a Leibniz algebra. 

\begin{Def}
	A Leibniz algebra morphism is a $\Bbbk-$linear map $f\colon L\rightarrow L'$ which is compatible with the bracket map, i.e. 
	\begin{equation*}
	f\left[x,y\right]=\left[f(x),f(y)\right]
	\end{equation*}
	for all $x,y\in L$.
\end{Def}

The category of Leibniz algebras consist of Leibniz algebras as objects and Leibniz algebra morphisms as morphisms. This category is denoted by $\Lbnz$.

\begin{Def}
	A Leibniz algebra with trivial bracket is called an Abelian (or singular) Leibniz algebra.
\end{Def}

\begin{Def}
	For any Leibniz algebras $L$ and $L'$ a Leibniz action of $L$ on $L'$ consist of two bilinear maps $\Lambda\colon L\times L' \rightarrow L', (x,a)\mapsto x\cdot m$ and $\rho\colon L'\times L \rightarrow L', (a,x)\mapsto m\cdot x$ satisfying
	\begin{enumerate}[label=\textbf{(\roman{*})}, leftmargin=2cm]
		\item	$x\cdot \left[ m,n \right]=\left[ x\cdot m,n \right]-\left[ x\cdot n,m \right]$,
		\item	$\left[ m,x\cdot n \right]=\left[ m\cdot x,n \right]-\left[ m,n \right]\cdot x$,
		\item	$\left[ m,n\cdot x \right]=\left[ m,n \right]\cdot x-\left[ m\cdot x,n \right]$,
		\item	$x\cdot \left( y\cdot m \right)=\left[ x,y \right]\cdot m-\left( x\cdot m \right)\cdot y$,
		\item	$x\cdot \left( m\cdot y \right)=\left( x\cdot m \right)\cdot y-\left[ x,y \right]\cdot m$,
		\item	$m\cdot \left[ x,y \right]=\left( m\cdot x \right)\cdot y-\left( m\cdot y \right)\cdot x$	
	\end{enumerate}
for all $x,y\in L$ and $m,n\in L'$.
\end{Def}

Let $L$ and $L'$ be two Leibniz algebras. A split extension of $L$ by $L'$ is a short exact sequence
\[\xymatrix{
	\mathcal{E} : & 0 \ar[r] & L' \ar@{->}[r]^{i} & E \ar@{->}[r]^{p} & L %\ar@/_/[l]_{s} 
	\ar[r] & 0 }\]
in $\Lbnz$ with a Leibniz algebra morphism $s\colon L\rightarrow E$ such that $ps=1_L$. Here, note that $p$ is surjective and $\ker p=i$. Given a split extension of $L$ by $L'$,  we get derived actions of $L$ on $L'$ defined by
\begin{align*}
x\cdot m=&\left[s(x),m\right]\\
m\cdot x=&\left[m,s(x)\right]
\end{align*}
for any $x\in L$ and $m\in L'$. Let a split extension 
\[\xymatrix{
	\mathcal{E} : & 0 \ar[r] & L' \ar@{->}[r]^{i} & E \ar@{->}[r]_{p} & L \ar@<-.7ex>[l]_{s} \ar[r] & 0 }\]
is given. Then by using the bijection 
\[\begin{array}{rcccl}
\theta & \colon & L'\times L & \longrightarrow & E \\
&        & \left(m,x\right)      & \longmapsto     & m+s(x)
\end{array}\]
we can define a Leibniz algebra structure on $L'\times L$ as follows:

\[\left[\left(m,x\right),\left(n,y\right)\right]=\left(\left[m,n\right]+m\cdot y+x\cdot n, \left[x,y\right]\right)\]
for all $x,y\in L$ and $m,n\in L'$. The inverse of the function $\theta$ is defined by

\[\begin{array}{rcccl}
\theta^{-1} & \colon & E & \longrightarrow & L'\times L \\
&        & e      & \longmapsto     & \theta^{-1}(e)=(e-sp(e),p(e))
\end{array}\]
for all $e\in E$. Thus $L'\times L$ cartesian product set becomes a Leibniz algebra which is called by semi-direct product of Leibniz algebras and denoted by $L'\rtimes L$.

For any Leibniz algebra $L$, the obvious action of $L$ on itself corresponds to the extension

\[\xymatrix{\mathcal{L} : & 0 \ar[r] & L \ar@{->}[r]^-{i} & L\rtimes L \ar@{->}[r]_-{p} & L \ar@<-.7ex>[l]_-{s} \ar[r] & 0 }\]

where $i(l)=(l,0)$, $p(l,l_1)=l_1$ and $s(l)=(0,l)$.

Now, we can give the definition of crossed modules of Leibniz algebras due to Porter \cite{Por}.

\begin{Def}\cite{Por}
Let $L_0$ and $L_1$ be two Leibniz algebras. Given a split extension 
\[\xymatrix{\mathcal{L} : & 0 \ar[r] & L_1 \ar@{->}[r]^-{i} & L_1\rtimes L_0 \ar@{->}[r]_-{p} & L_0 \ar@<-.7ex>[l]_-{s} \ar[r] & 0 }\]
of $L_0$ by $L_1$ and a Leibniz algebra morphism $\partial\colon L_1 \rightarrow L_0$, $\partial$ is called a \textbf{crossed module} if $(1_{L_{1}},\partial)$ and $(\partial,1_{L_{0}})$ are both split extension morphisms in $\Lbnz$.
\end{Def}
\[\xymatrix{
\mathcal{L}_1 : & 0 \ar[r] & L_1 \ar@{->}[r]^---{i} \ar[d]_{1_{L_{1}}} & L_1\rtimes L_1 \ar@{->}[r]^---{p} \ar[d]^{(1_{L_{1}},\partial)} & L_1  \ar[r] \ar[d]^{\partial} & 0 \\
\mathcal{L} : & 0 \ar[r] & L_1 \ar@{->}[r]^---{i} \ar[d]_{\partial} & L_1\rtimes L_0 \ar@{->}[r]^---{p} \ar[d]^{(\partial,1_{L_{0}})} & L_0 \ar[d]^{1_{L_{0}}} \ar[r] & 0 \\
\mathcal{L}_0 : & 0 \ar[r] & L_0 \ar@{->}[r]^---{i} & L_0\rtimes L_0 \ar@{->}[r]^---{p} & L_0  \ar[r] & 0 \\
}\]

A crossed module is denoted by $(L_1,L_0,\partial)$. It is more practical to have a description in terms of actions and Leibniz bracket. We recall the definitions from \cite{asa} and \cite{Por2}.

\begin{Prop}
	A crossed module of Leibniz algebras is a Leibniz algebra morphism $\partial\colon L_1 \rightarrow L_0$ with actions of $L_0$ on $L_1$ satisfying the following conditions for all $l_0\in L_0$ and $l_1,l'_1\in L_1$
	\begin{enumerate}[label=\textbf{(LXM\arabic{*})}, leftmargin=2cm]
		\item $\partial \left( {{l}_{0}}\cdot {{l}_{1}} \right)=\left[ {{l}_{0}},\partial \left( {{l}_{1}} \right) \right],\partial \left( {{l}_{1}}\cdot {{l}_{0}} \right)=\left[ \partial \left( {{l}_{1}} \right),{{l}_{0}} \right],$
		\item ${{l}_{1}}\cdot \partial \left( {{l}_{1}}' \right)=\left[ {{l}_{1}},{{l}_{1}}' \right],\partial \left( {{l}_{1}}' \right)\cdot {{l}_{1}}=\left[ {{l}_{1}}',{{l}_{1}} \right].$
	\end{enumerate}	
\end{Prop}

\begin{Prop}
	If $(L_1,L_0,\partial)$ is a crossed module, then $\ker \partial$ is an Abelian Leibniz algebra.
\end{Prop}
\begin{Prf}
It can easily be shown by using the crossed module condition \textbf{(LXM2)}.	
\end{Prf}

For any two crossed module $(L_1,L_0,\partial)$ and $(M_1,M_0,\delta)$ let $f_1\colon L_1 \rightarrow M_1$ and $f_0\colon L_0\rightarrow M_0$ be two Leibniz algebra morphisms. Then $(f_1,f_0)$ is called a crossed module morphism if the following conditions hold for all $l_0\in L_0$ and $l_1\in L_1$:
\begin{enumerate}[label=\textbf{(\roman{*})}, leftmargin=2cm]
	\item ${{f}_{0}}\circ \partial =\delta \circ {{f}_{1}}$,
	\item ${{f}_{1}}\left( {{l}_{0}}\cdot {{l}_{1}} \right)={{f}_{0}}\left( {{l}_{0}} \right)\cdot {{f}_{1}}\left( {{l}_{1}} \right)$,
	\item ${{f}_{1}}\left( {{l}_{1}}\cdot {{l}_{0}} \right)={{f}_{1}}\left( {{l}_{1}} \right)\cdot {{f}_{0}}\left( {{l}_{0}} \right)$
\end{enumerate}

Thus the category $\LX$ of Leibniz crossed modules can be constructed. The objects of this category are Leibniz crossed modules and morphisms are crossed module morphisms.

A groupoid is a category in which every morphism is an isomorphism. Let  $G$  be a groupoid. We write $\Ob(G)$ for the set of objects of $G$ and write $G$ for the set of morphisms. We also identify $\Ob(G)$ with the set of identities of $G$ and so an element of $\Ob(G)$ may be written as $x$  or  $1_x$  as convenient. We write $d_0, d_1 \colon G\rightarrow \Ob(G)$ for the source and target maps, and, as usual, write $G(x,y)$ for $d_0^{-1}(x)\cap d_1 ^{-1}(y)$, for $x,y\in \Ob(G)$. The composition  $h\circ g$ of two elements of $G$ is defined if and only if $d_0(h)=d_1(g)$, and so the map $(h,g)\mapsto h\circ g$ is defined on the pullback $G{_{d_0}\times_{d_1}} G$ of $d_0$  and $d_1$. The \emph{inverse} of $g\in G(x,y)$ is denoted by $g^{-1}\in G(y,x)$. If  $x\in \Ob(G) $, we write $\St_Gx$  for $d_0^{-1}(x) $  and  call the \emph{star} of $G$ at $x$.

A groupoid $G$ is \emph{transitive (resp. simply transitive, 1-transitive and totally intransitive)} if $G(x,y)\neq\emptyset$ (resp. $G(x,y)$ has no more than one element, $G(x,y)$ has exactly one element and $G(x,y)=\emptyset$) for all $x,y\in \Ob(G)$ such that $x\neq y$.

\section{Internal categories in $\Lbnz$}

\begin{Def}
	Let $\bmC$ be an arbitrary category with pullbacks. An internal category $C$ in $\bmC$ is a category in which the initial and final point maps $d_0,d_1\colon C\rightarrow \Ob(C)$, the object inclusion map $\epsilon\colon \Ob(C)\rightarrow C$ and the partial composition $\circ\colon C {_{d_{0}}}\times_{d_{1}} C\rightarrow C, (a,b)\mapsto a\circ b$ are the morphisms in the category $\bmC$.
\end{Def}

Let $G$ be an internal category in $\C$. If there exist a morphism $g'\in G$ such that $g\circ g'=\epsilon d_1 (c)$ and $g'\circ g=\epsilon d_0(c)$ for all morphisms $g\in G$, then $G$ is called an internal groupoid and $g'$ is called the inverse of $g$ which is denoted by $g^{-1}$.

Let $G$ be an internal category in the category $\Lbnz$ of Liebniz algebras. Then $G$ and $\Ob(G)$ are Leibniz algebras and the structural maps $(d_0,d_1,\epsilon,\circ)$ are Leibniz algebra morphisms. So we can give the following proposition.

\begin{Prop}
	Let $G$ be an internal category in $\Lbnz$. Then for all $x,y\in\Ob(G)$ and $g,g'\in G$
	\begin{enumerate}[label=\textbf{(\roman{*})}, leftmargin=2cm]
		\item ${{d}_{0}}\left( \left[ g,g' \right] \right)=\left[ {{d}_{0}}\left( g \right),{{d}_{0}}\left( g' \right) \right]$,
		\item ${{d}_{1}}\left( \left[ g,g' \right] \right)=\left[ {{d}_{1}}\left( g \right),{{d}_{1}}\left( g' \right) \right]$,
		\item $\varepsilon \left( \left[ x,y \right] \right)=\left[ \varepsilon \left( x \right),\varepsilon \left( y \right) \right]$, i.e. ${{1}_{\left[ x,y \right]}}=\left[ {{1}_{x}}{{,1}_{y}} \right]$,
		\item ${{\left[ g,g' \right]}^{-1}}=\left[ {{g}^{-1}},{{\left( g' \right)}^{-1}} \right]$.
	\end{enumerate}
\end{Prop}

Also note that the operation $\circ$ being a Leibniz algebra morphism implies that 
\begin{align*}
\left( h\circ g\right)+\left( h'\circ g'\right)&=\left( h+h' \right)\circ \left( g,g' \right)\\
\left[ h\circ g,h'\circ g' \right]&=\left[ h,h' \right]\circ \left[ g,g' \right]
\end{align*}
for all $g,g',h,h'\in G$ such that ${{d}_{1}}\left( g \right)={{d}_{0}}\left( h \right)$ and ${{d}_{1}}\left( g' \right)={{d}_{0}}\left( h' \right)$. These identities are called interchange laws. An application of the interchange laws is that the composition can be expressed by the addition as follows: for $g,h\in G$ such that ${{d}_{1}}\left( g \right)={{d}_{0}}\left( h \right)$

\[h\circ g = h -\varepsilon d_0(h)+g = g -\varepsilon d_0(h) + h. \]

Clearly, one can see that any internal category in $\Lbnz$ is an internal groupoid. Indeed, for any $g\in G$, $g^{-1}=\epsilon d_0(g)-g+\epsilon d_1(g)$ is the inverse morphism of $g$. Hence, we will use internal groupoid instead of internal category.

\begin{Exam}
	Every Abelian Leibniz algebra $L$ is an internal groupoid in $\Lbnz$ with algebra of objects $\Ob(L)$ is trivial, i.e. singleton.
\end{Exam}

\begin{Exam}
	Let $L$ be a Leibniz algebra. Then $L\times L$ becomes an internal groupoid over $L$ in $\Lbnz$. Here $d_0(l,l')=l$, $d_1(l,l')=l'$, $\epsilon(l)=(l,l)$ and the composition $(l',l'')\circ(l,l')=(l,l'')$ for all $l,l',l''\in L$.
\end{Exam}

\begin{Prop}
	Let $G$ be an internal groupoid in $\Lbnz$. Then $\St_G 0=\ker d_0$ is an ideal of $G$.
\end{Prop}
\begin{Prf}
	It can be shown by an easy calculation. 
\end{Prf}

\begin{Lem}\label{lemA}
Let $G$ be an internal groupoid in $\Lbnz$. If $g_1\in \ker d_0$ and $g_2\in \ker d_1$, then
\[\left[g_1,g_2\right]=\left[g_2,g_1\right]=0\]
\end{Lem}
\begin{Prf} Assume that $g_1\in \ker d_0$ and $g_2\in \ker d_1$. So compositions $g_1\circ \varepsilon(0)$ and $\varepsilon(0)\circ g_2$ are defined, where $\varepsilon(0)=0$ the identity element of addition operation and hence, of bracket operation. Then,
	\begin{align*}
\left[g_1,g_2\right]  & = \left[g_1\circ 0,0\circ g_2\right] \\
& = \left[g_1,0\right]\circ\left[0,g_2\right]  \tag{by interchange law}\\
& = 0\circ 0 \\
& = 0
\end{align*}	
\end{Prf}

\begin{Lem}\label{lemB}
Let $G$ be an internal groupoid in $\Lbnz$.	If $g_1\in \ker d_0$, then we have
\[ \left[g_1,\varepsilon d_1(g)\right]=\left[g_1,g\right] \]
and
\[ \left[\varepsilon d_1(g),g_1\right]=\left[g,g_1\right]. \]
\end{Lem}
\begin{Prf}
	Since $g_1\in \ker d_0$ and $g-\varepsilon d_1(g)\in \ker d_1$, one can prove the assertion of the Lemma by using Lemma \ref{lemA}.
\end{Prf}

Let $G$ and $H$ be two internal groupoids in $\Lbnz$. An internal groupoid morphism (internal functor) $f\colon G\rightarrow H$ is a morphism of underlying groupoids and Leibniz algebra morphism on both the algebra of morphisms and the algebra of objects. So, we can construct the category of internal groupoids in $\Lbnz$. This category may be denoted by $\Cat(\Lbnz)$ or $\Gpd(\Lbnz)$.

\begin{Theo}\label{equixmodint}
	The category $\LX$ of crossed modules in the category of Leibniz algebras and the category $\Cat(\Lbnz)$ of internal categories (groupoids) in the category of Leibniz algebras are naturally equivalent.
\end{Theo}

\begin{Prf}
	We sketch the proof and left to the reader some of calculations. Let $G$ be an internal groupoid in $\Lbnz$. Then $\ker d_0$ and $\Ob(G)$ are both Leibniz algebras and the restriction of the final point map 
	\[d_1 \colon \ker d_0 \rightarrow \Ob(G)\]
	is a Leibniz algebra morphism. Moreover $\Ob(G)$ acts on $\ker d_0$ by the maps
	\[\begin{array}{ccl}
	  \Ob(G)\times \ker d_0 & \longrightarrow & \ker d_0 \\
	 \left(x,g\right)      & \longmapsto     & x\cdot g = \left[\varepsilon(x),g\right]
	\end{array}\]
	and
	\[\begin{array}{ccl}
	\ker d_0 \times \Ob(G) & \longrightarrow & \ker d_0 \\
	\left(g,x\right)      & \longmapsto     & g\cdot x = \left[g,\varepsilon(x)\right]
	\end{array}\]
	These are derived actions, since these are obtained from the split extension
	\[\xymatrix@1{
		0 \ar[r] & \ker d_0 \ar@{->}[r]^{i} & G \ar@{->}[r]_-{d_{0}} & \Ob(G)  \ar@<-.7ex>[l]_-{\varepsilon} \ar[r] & 0 }\]
	Here we note that
	\begin{equation*}\label{semicong}
	\ker d_0 \rtimes \Ob(G)\cong G.
	\end{equation*}
	Also $\left(\ker d_0, \Ob(G), d_1\right)$ is a crossed module. Indeed,
	\begin{enumerate}[label=\textbf{(LXM\arabic{*})}, leftmargin=2cm]
		\item for all $x\in \Ob(G)$ and $g\in \ker d_0$
		\begin{align*}
		d_1(g\cdot x)    & = d_1\left(\left[g,\varepsilon(x)\right]\right) \\
		& = \left[d_1(g),d_1(\varepsilon(x))\right]   \tag{$d_1$ is a Leibniz algebra morphism}\\
		& = \left[d_1(g),x\right] 
		\end{align*}
		and similarly 
		\begin{align*}
		d_1(x\cdot g)& = d_1\left(\left[\varepsilon(x),g\right]\right) \\
		& = \left[d_1(\varepsilon(x)),d_1(g)\right]   \tag{$d_1$ is a Leibniz algebra morphism}\\
		& = \left[x,d_1(g)\right] 
		\end{align*}
		\item for all $g,g_1\in \ker d_0$
		\begin{align*}
		g\cdot d_1(g_1)  & = \left[g,\varepsilon(d_1(g_1))\right] \\
		& =  \left[g,g_1\right] \tag{by Lemma \ref{lemB}}
		\end{align*}
		and similarly
		\begin{align*}
		d_1(g_1)\cdot g   & = \left[\varepsilon(d_1(g_1)),g\right] \\
		& =  \left[g_1,g\right] \tag{by Lemma \ref{lemB}}
		\end{align*}
	\end{enumerate}
This construction defines a functor, $\eta$, from the category $\Cat(\Lbnz)$ of internal categories in the category of Leibniz algebras to the category $\LX$ of crossed modules in the category of Leibniz algebras.
\[ \eta \colon \Cat(\Lbnz) \longrightarrow \LX \]
Conversely, let $(L_1,L_0,\partial)$ be a crossed module of Leibniz algebras. Then $\left(L_1\rtimes L_0,L_0,d_0,d_1,\varepsilon,\circ\right)$ becomes an internal groupoid in $\Lbnz$, where $d_0(l_1,l_0)=l_0$, $d_1(l_1,l_0)=\partial(l_1)+l_0$, $\varepsilon(l_0)=(0,l_0)$, the composition
\[(l'_1,l'_0)\circ(l_1,l_0)=(l'_1+l_1,l_0)\]
for $l'_0=\partial(l_1)+l_0$ and the inverse $(l_1,l_0)^{-1}=(-l_1,\partial(l_1)+l_0)$. Now we need to show that these structural maps are Leibniz algebra morphisms. For all $(l_1,l_0),(l'_1,l'_0)\in L_1\rtimes L_0$
	\begin{align*}
	d_0\left(\left[(l_1,l_0),(l'_1,l'_0)\right]\right) & = d_0\left(\left(\left[l_1,l'_1\right]+l_1\cdot l'_0+l_0\cdot l'_1,\left[l_0,l'_0\right]\right)\right)  \\
	& =  \left[l_0,l'_0\right]  \\
	& = \left[d_0(l_1,l_0),d_0(l'_1,l'_0)\right] ,
	\end{align*} 
	\begin{align*}
	d_1\left(\left[(l_1,l_0),(l'_1,l'_0)\right]\right) & = d_1\left(\left(\left[l_1,l'_1\right]+l_1\cdot l'_0+l_0\cdot l'_1,\left[l_0,l'_0\right]\right)\right)  \\
	& = \partial\left(\left[l_1,l'_1\right]+l_1\cdot l'_0+l_0\cdot l'_1\right)+\left[l_0,l'_0\right]  \\
	& =  \partial\left[l_1,l'_1\right]+\partial(l_1\cdot l'_0)+\partial(l_0\cdot l'_1)+\left[l_0,l'_0\right] \\
	& = \left[\partial(l_1),\partial(l'_1)\right]+\left[\partial(l_1),l'_0\right]+\left[l_0,\partial(l'_1)\right]+\left[l_0,l'_0\right] \tag{by \textbf{(LXM1)}}\\
	& = \left[\partial(l_1),\partial(l'_1)+l'_0\right]+\left[l_0,\partial(l'_1)+l'_0\right] \tag{by bilinearity}\\
	& = \left[\partial(l_1),+l_0,\partial(l'_1)+l'_0\right] \tag{by bilinearity}\\	
	& = \left[d_1(l_1,l_0),d_1(l'_1,l'_0)\right]
	\end{align*}
	\begin{align*}
	\varepsilon\left(\left[l_0,l'_0\right]\right) & = \left(0,\left[l_0,l'_0\right]\right)  \\
	& =  \left(\left[0,0\right]+\left[0,l'_0\right]+\left[l_0,0\right],\left[l_0,l'_0\right]\right) \tag{ subs. $0=\left[0,0\right]+\left[0,l'_0\right]+\left[l_0,0\right]$ } \\
	& = \left[(0,l_0),(0,l'_0)\right]\\
	& = \left[\varepsilon(l_0),\varepsilon(l'_0)\right]
	\end{align*}
To see that the composition is a Leibniz algebra morphism, we need to verify the interchange law for bracket operation. Let $(l_1,l_0),(l'_1,l'_0),(l''_1,l''_0),(l'''_1,l'''_0)\in L_1\rtimes L_0$ such that $(l_1,l_0),(l'_1,l'_0)$ and $(l''_1,l''_0),(l'''_1,l'''_0)$ are composable, i.e. $l'_0=\partial(l_1)+l_0$ and $l'''_0=\partial(l''_1)+l''_0$. Then
\begingroup\makeatletter\def\f@size{9}\check@mathfonts
\def\maketag@@@#1{\hbox{\m@th\large\normalfont#1}}%
\begin{align*}
\left[(l'_1,l'_0)\circ(l_1,l_0),(l'''_1,l'''_0)\circ(l''_1,l''_0)\right] & = \left[(l'_1+l_1,l_0),(l'''_1+l''_1,l''_0)\right]\\
& = \left(\left[l'_1+l_1,l'''_1+l''_1\right]+\left(l'_1+l_1\right)\cdot l''_0+l_0\cdot\left(l'''_1+l''_1\right),\left[l_0,l''_0\right]\right)\\
& = \left(\left[l'_1,l'''_1\right]+\left[l_1,l'''_1\right]+\left[l'_1,l''_1\right]+\left[l_1,l''_1\right]+l'_1\cdot l''_0+l_1\cdot l''_0+l_0\cdot l'''_1+l_0\cdot l''_1,\left[l_0,l''_0\right]\right)\\
& = \left(\left[l'_1,l'''_1\right]+\partial(l_1)\cdot l'''_1 +l'_1\cdot\partial(l''_1)+\left[l_1,l''_1\right]+l'_1\cdot l''_0+l_1\cdot l''_0+l_0\cdot l'''_1+l_0\cdot l''_1,\left[l_0,l''_0\right]\right)\\
& = \left(\left[l'_1,l'''_1\right]+(\partial(l_1)+l_0)\cdot l'''_1 +l'_1\cdot(\partial(l''_1)+l''_0)+\left[l_1,l''_1\right]+l_1\cdot l''_0+l_0\cdot l''_1,\left[l_0,l''_0\right]\right)\\
& = \left(\left[l'_1,l'''_1\right]+l'_0\cdot l'''_1 +l'_1\cdot l'''_0+\left[l_1,l''_1\right]+l_1\cdot l''_0+l_0\cdot l''_1,\left[l_0,l''_0\right]\right)\\
& = \left(\left[l'_1,l'''_1\right]+l'_0\cdot l'''_1 +l'_1\cdot l'''_0,\left[l'_0,l'''_0\right]\right)\circ \left(\left[l_1,l''_1\right]+l_1\cdot l''_0+l_0\cdot l''_1,\left[l_0,l''_0\right]\right)\\
& = \left[(l'_1,l'_0),(l'''_1,l'''_0)\right]\circ\left[(l_1,l_0),(l''_1,l''_0)\right]
\end{align*} \endgroup
This shows that the composition $\circ$ is a morphism of Leibniz algebras. Thus $L_1\rtimes L_0$ becomes an internal groupoid on $L_0$ in $\Lbnz$. Above construction also defines a functor, $\delta$, from the category $\LX$ of crossed modules in the category of Leibniz algebras to the category $\Cat(\Lbnz)$ of internal categories in the category of Leibniz algebras.
\[ \delta \colon \LX \longrightarrow \Cat(\Lbnz) \]
It is straightforward to show that these functors, $\eta$ and $\delta$, gives a natural equivalence between the categories $ \LX $ and $ \Cat(\Lbnz) $, i.e. $ \eta\delta\simeq 1_{\LX} $ and $ \delta\eta\simeq 1_{\Cat(\Lbnz)} $.
\end{Prf}

\section{Coverings and actions of internal groupoids in $\Lbnz$}

First we will recall the definitions of coverings over groupoids from \cite{Br1}.

\begin{Def}(cf. \cite{Br1})
Let $p\colon\widetilde G\rightarrow G$ be a morphism of groupoids. Then $p$ is called a \emph{covering morphism} and $\widetilde{G}$  a \emph{covering groupoid} of $G$ if for each $\widetilde x\in \Ob(\widetilde G)$ the restriction  $\St_{\widetilde{G}}{\widetilde{x}} \rightarrow \St_{G}{p(\widetilde x)}$  is  bijective.	
\end{Def}

Assume that $p\colon \wtilde{G}\rightarrow G$ is a covering morphism. Then  we have a lifting function $S_{p}\colon G_{d_{0}}\times_{\Ob(p)}\Ob(\wtilde{G})\rightarrow \wtilde{G}$ assigning to the pair $(a,x)$ in the pullback $G_{d_{0}}\times_{\Ob(p)}\Ob(\wtilde{G})$ the unique element $b$ of $\St_{\widetilde{G}}\widetilde{x}$ such that $p(b)=a$. Clearly $S_{p}$ is inverse to $(p,d_{0})\colon  \widetilde{G}\rightarrow G_{d_{0}}\times_{\Ob(p)}\Ob(\widetilde{G})$. So it is stated that $p\colon  \widetilde{G}\rightarrow G$ is a covering morphism if and only if $(p,d_0)$ is a  bijection \cite{Br-Da-Ha}.

\begin{Def}
An internal groupoid morphism $p\colon \widetilde{G}\rightarrow G$ is a covering morphism if and only if $(p,d_0)\colon \widetilde{G}\rightarrow G_{d_{0}}\times_{\Ob(p)} \Ob(\widetilde{G})$ is an isomorphism in $\Lbnz$.
\end{Def}

A covering morphism $p\colon \widetilde{G}\rightarrow G$ is called \emph{transitive } if both $\widetilde{G}$ and $G$ are transitive.  A transitive covering morphism $p\colon\widetilde G\rightarrow G$ is called \emph{universal} if for every covering morphism $q\colon \widetilde{H}\rightarrow G$ there is a unique morphism of groupoids $\widetilde{p}\colon \widetilde G\rightarrow \widetilde{H}$ such that $q\widetilde{p}=p$ (and hence $\widetilde{p}$ is also a covering morphism), this is equivalent to that for $\widetilde{x}, \widetilde{y}\in \Ob({\widetilde G})$ the set $\widetilde{G}(\widetilde x, \widetilde y)$ has not more than one element.

\begin{Rem}\label{keriso}
Since for an internal groupoid $G$ in $\Lbnz$, the star $\St_G 0$ is also a Leibniz algebra, we have that if $p\colon \widetilde{G}\rightarrow G$ is a covering morphism of internal groupoids, then the restriction of $ p $ to the stars $\St_{\widetilde{G}}{\widetilde{0}} \rightarrow \St_{G}{0}$ is an isomorphism in $\Lbnz$.
\end{Rem} 

Let $p\colon \widetilde{G}\rightarrow G$ and $q\colon G'\rightarrow G$ be two coverings of $G$. A morphism $f\colon \widetilde{G}\rightarrow G'$ of coverings is a morphism of internal groupoids in $\Lbnz$ such that $qf=p$, i.e. following diagram is commutative.

\[\xymatrix{ \widetilde{G} \ar[rr]^{f} \ar[dr]_{p} &  & G' \ar[dl]^{q} \\ 
	& G &  }\]

Hence we can construct the category of covering internal groupoids of an internal groupoid $G$ in $\Lbnz$ which has covering morphisms of $G$ as objects and has morphisms of coverings as morphisms. This category will be denoted by $\LGdC/G$.

Recall that an action of a groupoid $G$ on a set $S$ via a function $\omega\colon S\rightarrow \Ob(G)$ is a function ${G}_{d_{0}}\times_\omega S\rightarrow S, (g,s)\mapsto g\bullet s$ satisfying the usual rules for an action: $\omega(g\bullet s)=d_1(g)$, $1_{\omega(s)}\bullet s=s$ and $(h\circ g)\bullet s=h\bullet (g\bullet s)$ whenever $h\circ g$ and $g\bullet s$ are defined. A morphism $f\colon (S,\omega)\rightarrow (S',\omega')$ of such actions is a function $f\colon S\rightarrow S'$ such that $w'f=w$ and $f(g\bullet s)=g\bullet f(s)$ whenever $g\bullet s$ is defined. This gives a category $\Act(G)$ of actions of $G$ on sets. For such an action, the action groupoid $G\ltimes S$ is defined to have object set $S$, morphisms the pairs $(g,s)$ such that $d_0(g)=\omega(s)$, source and target maps $d_0(g,s)=s$, $d_1(g,s)=g\bullet s$, and the composition \[(g',s')\circ (g,s)=(g\circ g',s)\] whenever $s'=g\bullet s$. The projection $q\colon G\ltimes S\rightarrow G, (g,s)\mapsto s$ is a covering morphism of groupoids and the functor assigning this covering morphism to an action gives an equivalence of the categories $\Act(G)$ and $\Cov/G$. Following equivalence of the categories was given in \cite{Br-Mu1}.

\begin{Prop} (cf. \cite{Br-Mu1}) The categories $\GGdC/G$ and $\GGdA(G)$ are equivalent.\end{Prop}

\begin{Def}
Let $G$ be an internal groupoid in $\Lbnz$. An action of the internal groupoid $G$ on a Leibniz algebra $L$ via $\omega$ consists of a Leibniz algebra morphism $\omega\colon L\rightarrow \Ob(G)$ from $L$ to the Leibniz algebra of objects $\Ob(G)$ and a Leibniz algebra morphism
\[\begin{array}{ccl}
{G} {_{d_{0}}\times_\omega} L & \longrightarrow & L \\
 \left(g,l\right)      & \longmapsto     & g\bullet l,
\end{array}\]	
which is called the action, satisfying
\begin{enumerate}[label=\textbf{(A\arabic{*})}, leftmargin=2cm]
	\item $\omega(g\bullet l)=d_1(g)$, 
	\item $1_{\omega(l)}\bullet l=l$,
	\item $(h\circ g)\bullet l=h\bullet (g\bullet l)$ ,
\end{enumerate}
whenever $h\circ g$ and $g\bullet l$ are defined.
\end{Def}

Note that the action being a Leibniz algebra morphism implies the following so called interchange laws:
\begin{equation*}
(g\bullet l)+(g'\bullet l')=(g+ g')\bullet(l+l') \label{interchangeaction}
\end{equation*}
\begin{equation*}
[(g\bullet l),(g'\bullet l')]=[g,g']\bullet[l,l'] \label{interchangeaction2}
\end{equation*}
for all $g,g'\in G$ and $l,l'\in L$, whenever both sides are defined. 

A morphism $f\colon (L,\omega)\rightarrow (L',\omega')$ of such actions is a morphism $f\colon L\rightarrow L'$ of Leibniz algebras such that $\omega' f=\omega$. This gives a category $\LGdA(G)$ of actions of $G$ on Leibniz algebras. 

For an action of $G$ on a Leibniz algebra $L$ via $\omega$, the action groupoid $G\ltimes L$ has a Leibniz algebra structure defined by \[(g,l)+(g',l')=(g+g',l+l'),\]
\[\left[(g,l),(g',l')\right]=\left(\left[g,g'\right],\left[l,l'\right]\right)\]
and with this operations $G\ltimes L$ becomes an internal groupoid in $\Lbnz$.

\begin{Prop}
Let $G$ be an internal groupoid in $\Lbnz$. The categories $\LGdC/G$ and $\LGdA(G)$ are equivalent.
\end{Prop}

\begin{Prf}
Let $p\colon \widetilde{G}\rightarrow G$ be a covering morphism in $\Cat(\Lbnz)$. Then $G$ acts on $\Ob(\widetilde{G})$ via $p_0\colon \Ob(\widetilde{G}) \rightarrow \Ob(G)$ by
\[\begin{array}{ccl}
{G} {_{d_{0}}\times_{p_{0}}} \Ob(\widetilde{G}) & \longrightarrow & \Ob(\widetilde{G}) \\
\left(g,\widetilde{x}\right)      & \longmapsto     & g\bullet \widetilde{x}=d_1(\widetilde{g}),
\end{array}\]
where $\widetilde{g}$ is the unique lifting of $g$ with initial point $\widetilde{x}$. It is easy to verify that this map is an action and a Leibniz algebra morphism, since $p$ is a Leibniz algebra morphism.

Conversely, let $G$ acts on a Leibniz algebra $L$ via $\omega\colon L\rightarrow \Ob(G)$. Then $q\colon G\ltimes L\rightarrow G, (g,l)\mapsto g$ is a covering morphism in $\Cat(\Lbnz)$. It is straightforward to confirm that these constructions defines the intended natural equivalence.
\end{Prf}

\begin{Exam}
Let $G$ be an internal groupoid in $\Lbnz$. Then $1_G\colon G\rightarrow G$ is a covering morphism in $\Cat(\Lbnz)$. The corresponding action to $1_G$ is constructed as follows: $G$ acts on $\Ob(G)$ via $1_{\Ob(G)}\colon \Ob(G)\rightarrow \Ob(G)$ where the action is 
\[\begin{array}{ccl}
{G} {_{d_{0}}\times_{1_{\Ob(G)}}} \Ob(G) & \longrightarrow & \Ob(G) \\
\left(g,x\right) & \longmapsto     & g\bullet x=d_1(g).
\end{array}\]
In this case the action groupoid 
\[G\ltimes \Ob(G)=\left\{\left(g,x\right) ~|~ d_0(g)=x\right\}\]
is isomorph to $G$ as an internal groupoid in $\Lbnz$, i.e., $G\ltimes \Ob(G) \cong G$.
\end{Exam}

\section{Covering crossed modules in $\Lbnz$}

The notion of coverings for crossed modules in the category of groups is introduced in \cite{Br-Mu1}. In a similar way, by using the equivalence of the categories $\Cat(\Lbnz)$ and $\XMod(\Lbnz)$, we can interpret the notion of coverings in $\XMod(\Lbnz)$.

\begin{Def}
Let $(L_1,L_0,\partial)$ and $(\widetilde{L_1},\widetilde{L_0},\widetilde{\partial})$ be two crossed modules of Leibniz algebras and $p_1\colon \widetilde{L_1}\rightarrow L_1$, $p_0\colon \widetilde{L_0}\rightarrow L_0$ be Liebniz algebra morphisms such that  $(p_1,p_0) \colon (\widetilde{L_1},\widetilde{L_0},\widetilde{\partial}) \rightarrow (L_1,L_0,\partial)$ is a crossed module morphism. If $p_1\colon \widetilde{L_1}\rightarrow L_1$ is an isomorphism of Leibniz algebras, then we say that $(\widetilde{L_1},\widetilde{L_0},\widetilde{\partial})$ is a covering crossed module of $(L_1,L_0,\partial)$ and that $(p_1,p_0)$ is a covering morphism of crossed modules.
\end{Def}  

\begin{Exam}
Let $(L_1,L_0,\partial)$ be a crossed modules of Leibniz algebras. Then $(1_{L_{1}},1_{L_{0}}) \colon (L_1,L_0,\partial) \rightarrow (L_1,L_0,\partial)$ is a covering.
\end{Exam}

Let $(p_1,p_0) \colon (\widetilde{L_1},\widetilde{L_0},\widetilde{\partial}) \rightarrow (L_1,L_0,\partial)$ and $(q_1,q_0) \colon (L'_1,L'_0,\partial') \rightarrow (L_1,L_0,\partial)$ be two coverings of $ (L_1,L_0,\partial) $. A morphism of coverings is a crossed module morphism $(f_1,f_0) \colon (\widetilde{L_1},\widetilde{L_0},\widetilde{\partial}) \rightarrow(L'_1,L'_0,\partial')$ such that $(q_1,q_0)\circ(f_1,f_0)=(p_1,p_0)$, i.e. $q_1 f_1=p_1$ and $q_0 f_0=p_0$. Now we can construct the category of coverings of $(L_1,L_0,\partial)$ which will be denoted by $\LXMC(L)$.

\begin{Prop}
Let $(L_1,L_0,\partial)$ be a crossed module of Leibniz algebras and $G$ be the corresponding internal groupoid	according to Theorem \ref{equixmodint}. Then the category $\LXMC(L)$ of coverings of $(L_1,L_0,\partial)$ and the category  $\LGdC/G$ covering internal groupoids of $G$ are equivalent.
\end{Prop}

\begin{Prf}
Let $p\colon \widetilde{G}\rightarrow G$ be a covering in $\Cat(\Lbnz)$ and $(\widetilde{L_1},\widetilde{L_0},\widetilde{\partial})$ be the corresponding crossed modules to $\widetilde{G}$. Then by Theorem \ref{equixmodint}, $\St_{\widetilde{G}}\widetilde{0}=\ker \widetilde{d_0}=\widetilde{L_1}$ and $\St_G 0=\ker d_0=L_1$. Since $p$ is a covering then by Remark \ref{keriso} the restriction of $p$ on $\widetilde{L_1}$ defines an isomorphism $\widetilde{L_1}\cong L_1$. Hence $(\widetilde{L_1},\widetilde{L_0},\widetilde{\partial})$ is a covering crossed module of $(L_1,L_0,\partial)$.

Conversely, let $(p_1,p_0)\colon(\widetilde{L_1},\widetilde{L_0},\widetilde{\partial})\rightarrow(L_1,L_0,\partial)$ be a covering of $(L_1,L_0,\partial)$ and $\widetilde{G}$ be the corresponding internal groupoid to $(\widetilde{L_1},\widetilde{L_0},\widetilde{\partial})$. Here $\widetilde{G}=\widetilde{L_1}\rtimes\widetilde{L_0}$, $\Ob(\widetilde{G})=\widetilde{L_0}$ and the corresponding internal groupoid morphism is $p=p_1\times p_0\colon\widetilde{G}\rightarrow G$. Let $x\in L_1$. Since $ \widetilde{L_1}\cong L_1 $ then there exist a unique $\widetilde{x}\in\widetilde{L_1}$ such that $p_1(\widetilde{x})=x$. Hence
\[\begin{array}{rcccl}
S_p &\colon& \left(L_1\rtimes L_0\right) {_{d_{0}}\times_{\Ob(p)}} \widetilde{L_0} & \longrightarrow & \widetilde{L_1}\rtimes\widetilde{L_0} \\
& & \left(\left(x,m\right),\widetilde{m}\right) & \longmapsto     & \left(\widetilde{x},\widetilde{m}\right)
\end{array}\]
defines an isomorphism of Leibniz algebras. One can easily see that these constructions are functorial and defines a natural equivalence between the categories $\LXMC(L)$ and $\LGdC/G$.
\end{Prf}

\section{Acknowledgement}
This work has been supported by Research Fund of the Aksaray University. Project Number:2016-019.

\small{

}

\end{document}